\documentclass{article}

\usepackage{arxiv}

\usepackage[utf8]{inputenc} 
\usepackage[T1]{fontenc}    
\usepackage{hyperref}       
\usepackage{url}            
\usepackage{booktabs}       
\usepackage{amsfonts}       
\usepackage{nicefrac}       
\usepackage{microtype}      
\usepackage{lipsum}
\usepackage{graphicx}
\graphicspath{ {./images/} }

\usepackage{amsthm,bbm,amssymb}
\usepackage{amsmath}

\usepackage[utf8]{inputenc}
\usepackage{amsmath,amssymb,natbib,xcolor,multicol,url,mhequ}
\usepackage{graphicx,bbm,xr}
\usepackage{booktabs,epstopdf,color}
\usepackage[space]{grffile}
\usepackage{lineno}
\usepackage[plain,noend]{algorithm2e}
\usepackage{multirow}
\DeclareMathOperator\supp{supp}

\DeclareMathOperator{\bbP}{\mathbb P}
\DeclareMathOperator{\E}{\mathbb E}
\DeclareMathOperator{\Var}{\mbox{Var}}

\def\T{\mathrm{\scriptscriptstyle{T}}}

\def\T{{\mathrm{\scriptscriptstyle T} }}

\newcommand{\be}{\begin{equs}}
\newcommand{\ee}{\end{equs}}

\numberwithin{equation}{section}
\theoremstyle{plain}
\newtheorem{theorem}{Theorem}

\newtheorem{corollary}{Corollary}
\newtheorem{lemma}{Lemma}

\usepackage{natbib}

\title{Adaptive posterior convergence in sparse high dimensional clipped generalized linear models}

\author{
 Biraj Subhra Guha \\
  Department of Statistics \\
  Texas A\&M University\\
  College Station, TX 77843 \\
  \texttt{birajsguha@stat.tamu.edu} \\
   \And
 Debdeep Pati \\
  Department of Statistics \\
  Texas A\&M University\\
  College Station, TX 77843 \\
  \texttt{debdeep@stat.tamu.edu} \\

}

\begin{document}
\maketitle
\title{Adaptive posterior convergence in sparse high dimensional clipped generalized linear models\thanksref{T1}}


\begin{abstract}
We develop a framework to study posterior contraction rates in sparse high dimensional generalized linear models (GLM).  We introduce a new family of GLMs, denoted by clipped GLM,  which subsumes many standard GLMs and makes minor modification of the rest.  With a sparsity inducing prior on the regression coefficients, we delineate sufficient conditions on true data generating density that leads to minimax optimal rates of posterior contraction of the coefficients in $\ell_1$ norm.   Our key contribution is to develop sufficient conditions commensurate with the geometry of the clipped GLM family, propose prior distributions which do not require any knowledge of the true parameters and avoid any assumption on the growth rate of the true coefficient vector.  
\end{abstract}

\section{Introduction}
The GLM \citep{mccullagh2019generalized} is a flexible generalization of ordinary linear regression that allows for response variables to accommodate error distributions which are non-additive and non-Gaussian. The GLM generalizes linear regression by allowing the linear model to be related to the response variable via a link function. Although primarily restricted to a lower dimensional setting,  Bayesian approaches for GLM  has been very popular from the 90's with the advent of Markov chain Monte Carlo \citep{dey2000generalized}.   

The emergence of more sophisticated data acquisition techniques in gene expression microarray, among many other fields, triggered the development of innovative statistical methods \citep{friedman2001elements,buhlmann2011statistics,hastie2015statistical} in the last decade, that help in analyzing large scale datasets.  The overarching  goal is to identify relevant predictors associated with a response out of a large number of predictors, but only with a smaller number of samples. This {\em large $p$, small $n$} paradigm is arguably the most researched topic in the last decade.  Primarily focusing on the linear models, statisticians have devised a number of penalized regression techniques for estimating $\beta$ in $p \gg n$ setting under the assumption of sparsity, with accompanying theoretical justification of optimal estimation, prediction and selection consistency; refer to \cite{tibshirani1996regression,fan2001variable,efron2004least,zou2005regularization,candes2007dantzig,zou2006adaptive,belloni2011square}.  Pioneering extensions of penalization based methods have been made for generalized linear models \citep{friedman2010regularization}, but existing results on theoretical guarantees for high dimensional GLMs are relatively few.  \cite{van2008high} studied the oracle rate of
the empirical risk minimizer with the lasso penalty in high dimensional GLMs. More recently, \cite{abramovich2016model} derived convergence rates with respect to the Kullback--Leibler risk with a wide class of penalizing functions, which can be translated into convergence rates relative to the $\ell_2$-norm under certain conditions.
 
 From a Bayesian standpoint, sparsity favoring mixture priors with separate control on the signal and noise coefficients have been proposed \citep{leamer1978regression,mitchell1988bayesian,george1995stochastic,george1997approaches,scott2010bayes,johnson2010use,narisetty2014bayesian,yang2016computational,rovckova2018spike}. Although in principle such methods can be used for generalized linear models, accompanying theoretical justification on optimal estimation in the high dimensional case is primarily available in the context of linear models \citep{castillo2012needles,castillo2015bayesian,gao2015general}.  

To the best of our knowledge, analogous results for generalized linear models in the high dimensional case are comparatively sparse, with the exception of \citet{jiang2007bayesian}.  However, special cases from the GLM family including high dimensional logistic regression using a pseudo likelihood \citep{atchade2017contraction} and high dimensional logistic regression using shrinkage priors \citep{wei2020contraction} are available.  \citet{jiang2007bayesian} operated in a high dimensional setting where the use of a Gaussian prior leads to a restrictive assumption on the growth of the true coefficients; refer to the assumptions of Theorem 1 in pg. 1493. \citet{atchade2017contraction} considered a Laplace-type prior for the coefficients which obviated the need for such a restriction, but their results are specific to logistic regression.

In this article, we develop a framework to study  posterior contraction in high dimensional clipped generalized linear models using complexity priors that involve a Laplace prior on the non-zero coefficients.   The clipped GLM class deviates slightly from the standard GLM construction in that we allow the effect of linear term $x ^{\T} \beta$ in the argument of the log-partition function to ``clip'' away from the singularities of the function.  Our clipped GLM directly subsumes high dimensional linear, polynomial and logistic regression, while also incorporating variants of Poisson, negative Binomial (and similar) regressions, which are identical from a practical standpoint to the standard Poisson/negative binomial regressions.

Our sufficient conditions are grouped into two categories: i) a set of identifiability and compatibility conditions based on the geometry of the clipped GLM, specified by the log-partition function that allows separation between models, and ii) an appropriate growth rate of scale parameter of the Laplace distribution that imposes appropriate penalty on the non-zero coefficients, along with an appropriate decay rate for the model weights that penalizes larger models. Existing literature \citep{jiang2007bayesian} on posterior contraction in GLMs requires growth rate assumptions on the true coefficient vector.  The crucial feature of our methodology is achieving adaptive, rate-optimal posterior contraction with respect to the data generation mechanism, while simultaneously avoiding any growth assumptions on the true coefficients.  

While our article was in final stages, we came across a dissertation by Seonghyun Jeong at NC State University under the supervision of Prof. Subhashis Ghosal, which considers posterior contraction in GLMs
using complexity priors on the model space in Chapter 4.  Their results make use of the same identifiability and compatibility assumptions as in \citep{castillo2015bayesian} to deliver optimal posterior contraction rates, albeit with a growth restriction on the true coefficient vector.  On the other hand, we do not require any growth assumption on the true coefficient vector.  
 Our assumptions  for obtaining adaptive rate-optimal posterior contraction are specifically designed for the clipped GLMs which can be viewed as appropriate generalization of the identifiability and compatibility assumptions of \cite{castillo2015bayesian} in the linear model case. Finally, the prior dependence on the 
true parameter can be completely eliminated making our results rate-adaptive.  

The remaining of the article is organized as follows. Section \ref{sec:GLM} introduces the construction of the clipped GLM family. Section \ref{sec:prior} details the sparsity favoring prior construction while section \ref{sec:gendata} entails the identifiability and compatibility assumptions on the data generating process and the choice of hyperparameters.  Section \ref{sec:rate} states our main results on adaptive rate-optimal posterior contraction. This is divided into three parts: a lower bound on the marginal likelihood, a result on controlling the effective sparsity of the posterior distribution and finally a truth-adaptive contraction rate theorem. The proofs are deferred to the Appendices \ref{sec:thm1}-\ref{sec:thm3} with the auxiliary results in Appendix \ref{sec:aux}. 

\subsection{Notations}
For reals $\zeta_1, \zeta_2$,  $\zeta_1 \precsim \zeta_1$ denotes $\zeta_1 \leq C_1 \zeta_2$ for an absolute constant $C_1$. Similarly, we define $\zeta_1 \gtrsim \zeta_2$. For sequences of real numbers $\left\{ \zeta_{1, n} \right\}$ and $\left\{ \zeta_{2, n} \right\}$, we say $\zeta_{1, n} = o(\zeta_{2, n})$ if $\zeta_{1, n}/\zeta_{2, n} \to 0$, and $\zeta_{1, n} = O(\zeta_{2, n})$ if we have $0 < C_2 \leq \underset{n \to \infty}{\liminf} \left( \zeta_{1, n}/\zeta_{2, n} \right) \leq \underset{n \to \infty}{\limsup} \left( \zeta_{1, n}/\zeta_{2, n} \right) \leq C_3$ for absolute constants $C_2, C_3$.

\section{Construction of the GLM family} \label{sec:GLM}
For both univariate and multivariate observations, one of the most widely used and well-structured family of models is the exponential family. One can refer to \cite{koopman1936distributions} and \cite{pitman1936sufficient} for the initial works on exponential families. We discuss this briefly with the example of univariate observations and real valued parameter. The exponential family takes the form
\begin{equation}\label{2_1}
    f(y\mid \theta) = h(y) \exp \left[ \theta T(y) - A(\theta) \right], \; y \in \mathcal{Y} \subset \mathbb{R},
\end{equation}
where $\theta \in \Theta \subset \mathbb{R}$ is the parameter of interest, $h(\cdot):\mathcal{Y} \to \mathbb{R}$ is called the base measure, $A(\cdot): \Theta \to \mathbb{R}$ is the convex log-partition function and $T(\cdot):\mathcal{Y} \to \mathbb{R}$ is called the sufficient statistic for estimating parameter $\theta$. This form is known as the canonical form of an exponential family. Many standard distributions, like the Bernoulli and Gaussian with known variance, Poisson, negative Binomial, among many others, follow model \eqref{2_1}. It is well known that the mean and variance of the sufficient statistic is given in terms of $A(\cdot)$, namely $\E[T(Y)] = A'(\theta), \Var(T(Y)) = A''(\theta)$. $A'(\cdot)$ and $A''(\cdot)$ are thus known as the mean and variance functions respectively, and 
$A'(\cdot)$ can be assumed to strictly increasing on its domain. An interesting property of exponential families is that it affords a neat expression of Kullback--Leibler(KL) divergence in terms of the Bregman divergence of log-partition function $A(\cdot)$. The Bregman divergence of convex function $A(\cdot)$ at $\theta_0$ from $\theta$ is given by $A(\theta) - A(\theta_0) - \left( \theta - \theta_0 \right) A'(\theta_0)$, which turns out to be the same expression for KL divergence of $\theta_0$ from $\theta$, which we denote by $\mathcal{D}(\theta_0| |\theta)$. These properties play a major role in dealing with exponential family distributions.  

A generalized linear model (GLM) assumes that the observation comes from an exponential family member as above, and models a function of the mean through a linear function of a covariate, i.e. as $x^{\T} \beta$, where $x$ represents a covariate and $\beta$ is the new parameter vector of interest. The said function, denoted by $g(\cdot): \mathrm{range}[A'(\cdot)] \to \mathbb{R}$, is termed as the link function. With $n$ data points and $d_n$ covariates, $X:n \times d_n$ makes up the design matrix, whose $i$-th row is denoted by $x_i^{\T}$. Thus, for every $i = 1, \dots n$, GLM prescribes the transition $\theta$ to $\beta$ as
\begin{equation}\label{2_2}
    g^{-1} \left(x_i^{\T} \beta \right) = A'(\theta), \,\mbox{equivalently}\, \theta = \left( g \circ A' \right)^{-1} \left( x_i^{\T} \beta \right).
\end{equation}
It is clear from the right hand side of \eqref{2_2} that GLM actually models the original parameter of the exponential family, but it does so indirectly, through the link function and $A'(\cdot)$. As we shall see in our next section \ref{sec:cGLM}, \eqref{2_2} motivates modeling the original parameter $\theta$ using $A''(\cdot)$, and not through $A'(\cdot)$, leading to the definition of clipping function $\eta(\cdot)$ and clipped GLM family.

\subsection{Introduction to clipped GLM}\label{sec:cGLM}
\begin{figure}[htp]
    \centering
    \includegraphics[width=8cm]{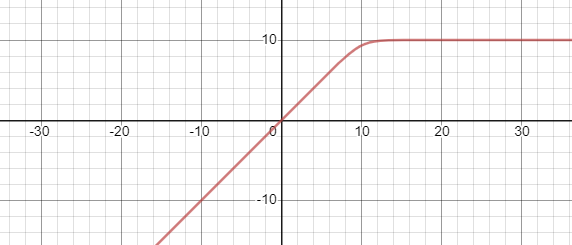}
    \caption{Graph of $y = 10 - \log \left( 1 + \exp ( - x + 10) \right)$}
    \label{fig:clipping}
\end{figure}

We now discuss in detail the clipped Generalized Linear Model (cGLM), which includes, but are not limited to, the distributions like Bernoulli, binomial with known number of trials, Poisson, negative binomial with known number of failures, exponential, Pareto with known minimum, Weibull with known shape, Laplace with known mean, chi-squared and Gaussian with known known variance. We start with the canonical rank-one exponential family of distributions, where the canonical parameter $\theta$ is expressed through a function of covariates. However, in contrast to GLM, we choose to represent
$$\theta = \eta \left( x^{\T} \beta \right),$$
where $\eta(\cdot)$, termed as the clipping function, depends only on $A''(\cdot)$. In cGLM, we consider log-partition functions $A(\cdot)$ that satisfy
\begin{itemize}
    \item $A''(\cdot)$ exists everywhere in the domain of $A(\cdot)$,
    \item $\mathcal{I}_A (b) := \{t \in \mathbb{R} : 0 \leq A''(t) \leq b \}$ is an interval on the real line for any $b \in (0, \infty]$.
\end{itemize}
All the standard examples of exponential families we discuss satisfy these simple properties. We now turn to the clipping functions we use in cGLM, which play an intermediary role, sitting between $A(\cdot)$ and the $i$-th linear term $x_i^{\T} \beta$. We motivate the choice of clipping functions by describing some examples. Since we work with $\beta \in \mathbb{R}^{d_n}$, the linear term $x_i^{\T} \beta$ belongs to $\mathbb{R}$, whereas the log-partition function $A$ can have strict interval subsets of the real line as their support. These types of log-partition functions have a single {\em pole} ($r_0$ such that $\lim_{x \to r_0} A(x) =\infty$) on the real line. Examples include:
\begin{itemize}
    \item \textit{Negative Binomial:} $A(t) = -q \log \left( 1 - \exp(t) \right), t < 0$ with $q$ denoting known number of failures. This shows $r_0 = 0$.
    \item \textit{Exponential:} $A(t) = - \log (-t), t < 0$ so that $r_0 = 0$.
    \item \textit{Pareto:} $A(t) = - \log(- 1 - t) + (1 + t) \log q_{\mathrm{min}}, t < -1$ with $q_{\mathrm{min}}$ denoting known minimum value. This shows $r_0 = -1$.
    \item \textit{Laplace:} $A(t) = - \log (- t/2), t < 0$ so that $r_0 = 0$. Mean is assumed to be zero.
\end{itemize}
Distributions like Bernoulli (or Binomial with known number of trials), Poisson and Gaussian (with known variance) have log-partition functions with entire real line as support. The clipping function's first role is to ensure that $\eta_i \equiv \eta \left( x_i^{\T} \beta \right)$, which acts as an argument to $A(\cdot)$ to have the same range as the domain of $A(\cdot)$. Its second role, which turns out to be the central point of our hyper-parameter assumption, is to control the growth of $A''(\cdot)$, specifically to allow a local quadratic majorizability of $A(\cdot)$. Bernoulli and Gaussian (with known variance) already enjoy the special status of having a universal bound on $A''(\cdot)$. Hence, for Poisson, which has $A(t) = \exp(t), t \in \mathbb{R}$, as well as the distributions that have a pole in their log-partition function, $\eta(\cdot)$ should be assumed to be playing the role of clipping the linear term $x_i^{\T} \beta$ away from $+ \infty$ and $r_0$ respectively, or $\pm \infty$ and poles in general cGLM members. We illustrate one possible set of choices of clipping function $\eta(\cdot)$ in the following examples. It is important to note their connection to the popular regression settings, which we shall delve into in \eqref{sec:reg}. 
\begin{itemize}
    \item \textit{Bernoulli:} $\eta(t) = t$, due to universal bound on $A''(\cdot)$.
    \item \textit{Negative binomial with known number of failures:} $\eta(t) = - \delta - \log \left( 1 + \exp ( - t - \delta) \right)$, where $\delta$ is a small positive absolute constant.
    \item \textit{Poisson:} $\eta(t) = \mathcal{C}_0 - \log \left( 1 + \exp ( - t + \mathcal{C}_0) \right)$, where $\mathcal{C}_0$ is large positive absolute constant (see figure \eqref{fig:clipping}, where $\mathcal{C}_0 = 10$).
    \item \textit{Exponential:} $\eta(t) = - \delta - \log \left( 1 + \exp ( - t - \delta) \right)$, where $\delta$ is a small positive absolute constant.
    \item \textit{Gaussian with known variance:} $\eta(t) = t$, due to universal bound on $A''(\cdot)$.
    \item \textit{Pareto with known minimum value:} $\eta(t) = - (1 + \delta) - \log \left( 1 + \exp ( - 1 - t - \delta) \right)$, where $\delta$ is a small positive absolute constant.
    \item \textit{Laplace with known mean:} $\eta(t) = - \delta - \log \left( 1 + \exp ( - t - \delta) \right)$, where $\delta$ is a small positive absolute constant.
\end{itemize}

Two points are crucial to note here. The clipping function $\eta(\cdot)$ can be defined as injective and Lipschitz, as all our examples show. These two properties play an important role in identifiability of the model, as is discussed in the next section. Secondly, the constants $\mathcal{C}_0, \delta$ are absolute, meaning that the practitioner should choose them before-hand, and their choice is totally independent of the observed data or the true value of the parameter in question. An example of such a choice would be $\delta = 10^{-4}$ and $\mathcal{C}_0 = 10^3$. We now summarize the defining properties of clipping functions $\eta(\cdot)$ used in cGLM, and their connection to $A(\cdot)$ through a the following simple and mild condition: \newline

\noindent \textbf{Clipping function condition:} There exists constant $\mathcal{M}_0(A) > 0$ depending on $A(\cdot)$, so that $\eta(\cdot)$ satisfies
\begin{equation}\label{2_3}
    \eta(\cdot) : \mathbb{R} \to \mathcal{I}_A \left( \frac{\mathcal{M}_0^2(A)}{2} \right), \mathrm{\, Lipschitz, \, injective}.
\end{equation}
\newline

We now describe our data-generating model. For $i = 1, \dots n$, $y_i \in \mathcal{Y} \subset \mathbb{R}$ are independent data points with $x_i \in \mathbb{R}^{d_n}$ as the covariate, $\beta \in \mathbb{R}^{d_n}$ as the parameter of interest and $\beta^{*}$ denoting the true parameter value. Let $\eta_i \equiv \eta \left( x_i^{\T} \beta \right), \eta_i^{*} \equiv \eta \left( x_i^{\T} \beta^{*} \right)$ and let $X$ denote the covariate matrix or design matrix, with the vector $x_i^{\T}$ representing the $i$-th row of $X$. The sufficient statistic is $T_i \equiv T \left( y_i \right)$, the base measure by $h(y_i)$ and the density for the $i$-th data point is denoted by $f\left(y_i \mid \eta_i \right)$. The $i$-th log-partition function  is denoted by $A\left( \eta_i \right)$. We denote by $S^{*}$ the true model, the non-zero co-ordinates of $\beta^{*}$. Also, we shall denote by $\supp{(\beta)}$ the set of non-zero entries in $\beta$, and by $\beta_S$ the same vector as $\beta$ with the co-ordinates in $S^c$ set to zero. $L_n(\eta, \eta^{*})$ stands for the log-likelihood ratio, which is expressed in terms of its two parts; $Z_n(\eta, \eta^{*})$ is the centered stochastic term and while $\mathcal{D}_n (\eta^{*}| |\eta)$ denotes the Kullback--Leibler(KL) divergence, both based on $y^{(n)}$. Thus we have the following:
\begin{equation}\label{2_4}
    \begin{split}
        f\left(y_i \mid \eta_i \right) &= h \left(y_i \right) \exp \left( T_i \eta_i - A \left( \eta_i \right) \right), \; i=1, \dots n, \\
         \mathcal{D}_i (\eta_i^{*}| |\eta_i) &:=  A\left(\eta_i \right) - A\left(\eta_i^{*} \right) -  \left( \eta_i - \eta_i^{*} \right) A' \left( \eta_i^{*} \right), \; \mathcal{D}_n (\eta^{*}| |\eta) := \sum_{i = 1}^n \mathcal{D}_i (\eta_i^{*}| |\eta_i), \\
        Z_i (\eta_i, \eta_i^{*}) &:= (T_i - \E T_i) \left(\eta_i - \eta_i^{*} \right), \; Z_n (\eta, \eta^{*}) := \sum_{i = 1}^n Z_i (\eta_i, \eta_i^{*}), \\
        L_n(\eta, \eta^{*}) &:= Z_n(\eta, \eta^{*}) - \mathcal{D}_n (\eta^{*}| |\eta).
    \end{split}
\end{equation}

\subsection{Connection of cGLM to regression settings}\label{sec:reg}
We briefly discuss how cGLM incorporates more commonly used high dimensional linear and non-linear regression settings. As we shall see, GLM and cGLM are interchangeable from the standpoint of practical implementation. Recall from \eqref{sec:cGLM} that we model the canonical parameter $\theta$ of the exponential family underlying cGLM as $\theta = \eta \left( x_i^{\T} \beta \right)$.
\begin{itemize}
    \item \textit{Linear regression with Gaussian error:} Since the native parameter, which is the mean, is the same as the canonical parameter for Gaussian, the choice of normal distribution with known variance for the exponential family in cGLM, alongside the valid choice of $\eta(t) = t, \; t \in \mathbb{R}$ as the clipping function, leads us to classical high dimensional linear regression (large $d_n$ and small $n$) with Gaussian errors. Here, $\mathcal{Y} = \mathbb{R}$.
    
    \item \textit{Logistic regression:} The native parameter here is probability of success $p \in (0, 1)$, while the canonical parameter is $\theta = \log(p/(1-p)) \in \mathbb{R}$. Thus, choosing Bernoulli for the exponential family and then, similar to linear regression, taking $\eta(t) = t, \; t \in \mathbb{R}$ as the clipping function, gives us the standard logistic regression setup. Here, $\mathcal{Y} = \{0, 1\}$.
    
    \item \textit{Poisson regression:} Denoting the native parameter in Poisson as $\nu > 0$, we see that the canonical parameter takes the form $\theta = \log \nu \in \mathbb{R}$. Standard Poisson regression would demand of us an identity clipping function alongside the choice of Poisson for the exponential family. However, because of \eqref{2_3}, we can allow $\eta(t) = \mathcal{C}_0 - \log \left( 1 + \exp ( - t + \mathcal{C}_0) \right), \; t \in \mathbb{R}$ for a large $\mathcal{C}_0 > 0$ of the practitioner's choosing. Refer to \eqref{fig:clipping} for a graph of this clipping function when $\mathcal{C}_0 = 10$. As is clear, we are allowing $\eta \left( x_i^{\T} \beta \right) $ to be approximately $ x_i^{\T} \beta$ i.e. linear, on $t < \mathcal{C}_0$, which is desired in Poisson regression, but clipping it to almost the  constant value $\mathcal{C}_0$ on $t \geq \mathcal{C}_0$. Intuitively for Poisson regression, $A(t) = \exp(t), t \in \mathbb{R}$ is already very large for moderately large $t$, hence allowing $t = \eta \left( x_i^{\T} \beta \right)$ to be large for large $\| \beta \|_1$ serves no extra purpose from a modelling perspective. Our choice of $\eta(\cdot)$ reflects this, maintaining negligible difference of GLM and cGLM from implementation perspective. Here, $\mathcal{Y} = \{0, 1, \dots \}$.
    
    \item \textit{Negative binomial regression with known number of failures:} The native parameter here is probability of success $p \in (0, 1)$, while the canonical parameter is $\theta = \log p \in (- \infty, 0)$. In contrast to the regression setups described above, standard negative binomial regression would require of us the clipping function $\eta(t) = - \log \left( 1 + q.\exp(t) \right)$, where $q \geq 1$ is the known number of failures, alongside choosing negative binomial for the exponential family. However, such a choice is unwarranted owing to \eqref{2_3}. Instead, we can go with $\eta(t) = - \delta - \log \left( 1 + \exp ( - t - \delta) \right), \; t \in \mathbb{R}$, $\delta > 0$ being a pre-fixed, small constant the practitioner decides upon. Our cGLM based choice of $\eta(\cdot)$, which almost completely mimics the GLM dictated choice, appropriately reflects the presence of pole at $r_0 = 0$ for negative binomial's log-partition function $A(t) = -q \log \left( 1 - \exp(t) \right), t < 0$. Here, $\mathcal{Y} = \{0, 1, \dots \}$.
\end{itemize}

\section{Construction of sparsity favoring prior}\label{sec:prior}
The sparsity favoring prior on the high dimensional $\beta$ is motivated by \cite{castillo2012needles,castillo2015bayesian} and follows the construction of spike-and-slab prior proposed in the early references \citep{leamer1978regression,mitchell1988bayesian,george1995stochastic,george1997approaches}. The crucial difference is in the slab part; we use a Laplace prior as in \cite{castillo2012needles,castillo2015bayesian}  instead of the more commonly used Gaussian slab. More recently, \cite{,johnson2010use} advocated the use of spike-and-non-local prior which has a better subset selection property compared to the spike-and-slab priors. However, our primary focus is in consistent estimation of $\beta$ and a spike-and-Laplace suffices in achieving this goal.  

The prior on parameter $\beta$ is induced through a prior on the duo $(S, \beta)$, where $S$ denotes a subset of $\{ 1, \dots d_n \}$. First, the prior on the dimension $0 \leq s \leq d_n$ is chosen to be $\omega_n(s) = C_n d_n^{-a_n s}, \; s = 0, \ldots, d_n$ with hyper-parameter $a_n > 0$, where $C_n$ is chosen to normalize the distribution. For any $\beta$ and  $S$ mentioned above, recall that  $\beta_S$ denotes the same vector $\beta$, but co-ordinates in $S^c$ set to $0$. With hyper-parameter $\lambda_n > 0$, the full prior is taken to be of the form
\begin{equation}\label{3_1}
   \begin{split}
       \Pi_n \left( S, \beta \right)  &:= \omega_n(|S|). {d_n \choose |S|}^{-1} . \; \left( \frac{\lambda_n}{2} \right)^{|S|}. \; \exp (-\lambda_n \|\beta_S\|_1) . \; \delta_0 \left( \beta_{S^{c}} \right) \\
       &= C_n. {d_n \choose |S|}^{-1} . \; \left( \frac{\lambda_n}{2 d_n^{a_n}} \right)^{|S|}. \; \exp (-\lambda_n \|\beta_S\|_1) . \; \delta_0 \left( \beta_{S^{c}} \right),
   \end{split}
\end{equation}
where $\|.\|_1$ denotes $\ell_1$-norm of Euclidean vectors, $|S|$ denotes cardinality of the set $S$ and $\delta_0$ denotes the degenerate distribution. The prior on the main parameter of interest, $\beta$, is given by
$$\Pi_n(\beta) := \sum_{S \subset \{ 1, \dots n \}} \Pi_n \left( S, \beta \right),$$
and the posterior probability of a general $B \subset \mathbb{R}^{d_n}$ is
$$\Pi_n(B \mid Y^{(n)}) := \frac{\int_{B} \exp [L_n(\eta, \eta^{*})] \Pi_n(\beta) d \beta}{\int \exp [L_n(\eta, \eta^{*})] \Pi_n(\beta) d \beta}.$$

\subsection{Prior on model size and the non-zero coefficients}
Our choice for model weights $\omega_n(\cdot)$ is special case of what is known as a complexity prior in \cite{castillo2015bayesian}. The prior is designed to down-weight models based on their larger sizes, and weight decrease is geometric in model dimension. We thus induce sparsity in the posterior through our prior choice. We point out that there are multiple ways of specifying and generalizing the prior we have used, specifically as in \cite{castillo2012needles} and \cite{castillo2015bayesian}, and they all share the central theme of exponential down-weighting of bigger models, and have the same effect on the posterior as our prior. We place independent Laplace signals for the non-zero coordinates. One can find dependent priors in the literature in this setup, for example in \cite{castillo2012needles}, but we choose to work with independent signals aiming to make our analysis neater. 

\section{Assumptions on data generating distribution and prior}\label{sec:gendata}
Our assumptions on the likelihood stem from that on the KL divergence term, while assumptions about the prior come from assumptions on the hyper-parameters $\lambda_n$ and $a_n$. These assumptions also dictate the possible values of true $\beta^{*}$, uniformly over which we shall state our results. In the first subsection, we present identifiability and compatibility (IC) conditions, and connect them to uniformly adaptive statements about the posterior. The second subsection is concerned with the choice of hyper-parameters that avoid any dependence of the prior on true $\beta^{*}$. 
We start by describing some order conditions, which shall help us define the rest of the assumptions.

\subsection{Order assumptions on sample size and parameter dimension:}
Since we work with a high dimensional problem, a natural condition is $d_n > n$ where $n, d_n \to \infty$. Now define a deterministic sequence of positive reals $\{ b_n \}$, such that
\begin{equation}\label{4_1}
     b_n = o \left( \frac{n}{\log d_n} \right).
\end{equation}
We shall focus on those true $\beta^{*}$'s whose sparsity $s_n^*$ satisfies $1 \leq s_n^{*} \leq b_n$. This gives us, among other things, the important relation: $(s_n^{*} \log d_n)/n \to 0$ as $n \to \infty$. It is also important that $s_n^{*} \nrightarrow 0$, which first gives us $s_n^{*} \log d_n \to \infty$, and second, forces us to have $\log d_n = o(n)$. This shows that $b_n = O(1)$ is a valid choice, satisfying \eqref{4_1}. We work with $n$ large enough so that $b_n \log d_n < n$ for all our calculations. Also, note that $d_n > n$ implies $3 b_n < d_n$ for large enough $n$.

\subsection{Identifiability and compatibility assumptions}
The ability of the log-likelihood term $L_n(\eta, \eta^{*})$ to create a separation between the true value of $\beta^{*}$ from any other $\beta$ is a fundamental criterion in posterior contraction analysis, and is termed as the identifiability criterion. Again, the natural measure of discrepancy in cGLM model is the Kullback--Leibler(KL) divergence $\mathcal{D}_n (\eta^{*}| |\eta)$, and since we work with Laplace signals in our prior, it is a natural demand to connect $\mathcal{D}_n (\eta^{*}| |\eta)$ with the $\ell_1$ distance, making them compatible. The requirements of compatibility and identifiability are simultaneously met by enforcing a lower bound on the KL divergence in terms of $\ell_1$ distance between the $\beta$'s i.e. $\| \beta_2 - \beta_1 \|_1$ for $\beta_1, \beta_2 \in \mathbb{R}^{d_n}$. We express this through the IC (Model) and IC (Dimension) assumptions, essentially requiring existence of a model $S$ and a dimension $s$, where $S \subset \{1, \dots d_n\}$ and $ s = 3 b_n, \dots d_n$ and they satisfy a certain lower bound property through the KL term. These assumptions not only generalize the compatibility assumptions made in \cite{castillo2015bayesian}, but also link them to identifiability of the truth. \newline

\noindent \textbf{IC (Model) Assumption:} There exists at least one non-null model $S \subset \{ 1, \dots d_n \}$ and the corresponding quantity $\phi_1(A, X, S) > 0$, such that for any $\beta_1, \beta_2 \in \mathbb{R}^{d_n}$, we have
$$\beta_{1S} \neq \beta_{2S}, \; \| \beta_{2S^c} - \beta_{1S^c} \|_1 < 7 \| \beta_{2S} - \beta_{1S} \|_1 \quad \Rightarrow \quad \mathcal{D}_n(\eta_1| |\eta_2) \geq \frac{n \phi_1^2(A, X, S)}{|S|} \| \beta_{2S} - \beta_{1S} \|_1^2.$$
The suffix $1$ of $\phi_1$ emphasizes we are working with constraints in the $\ell_1$ distance, as seen above. Intuitively, a general $\beta$, that is close to the truth $\beta^{*}$ in $\ell_1$ norm, will tend to have smaller absolute values in the true noise co-ordinates $S^{*c}$, and hence such a $\beta$ will tend to satisfy $\| \beta_{S^{*c}} - \beta^{*}_{S^{*c}} \|_1 < 7 \| \beta_{S^{*}} - \beta^{*}_{S^{*}} \|_1$ or equivalently $\| \beta_{S^{*c}} \|_1 < 7 \| \beta_{S^{*}} - \beta^{*} \|_1$. It is precisely in this scenario we shall need the IC (Model) assumption, i.e., $\phi_1(A, X, S^{*}) > 0$ so that the KL term creates a separation of the true and non-true $\beta$'s that are close in $\ell_1$ distance.  The IC (Model) assumption will be crucially used in our proof of Theorem \ref{thm:2}. \newline
 
Now consider the following subset of the parameter space:
$$\mathcal{B}_{1, n} := \left\{\beta \in \mathbb{R}^{d_n} : \phi_1 \left(A, X, \supp{(\beta)} \right) > 0 \right\}.$$
Based on the previous discussion, we would need true $\beta^{*} \in \mathcal{B}_{1, n}$, and due to IC(Model) assumption, $\mathcal{B}_{1, n}$ is non-null. Also, given any $A$ and $X$, the quantity $\phi_1 (A, X, S)$ can only take finitely many values as $S$ varies over subsets of $\{ 1, \dots d_n \}$, all of those values being positive for $S = S^{*}$. This gives us the quantity, for any non-null $\mathcal{B} \subset \mathcal{B}_{1, n}$,
\begin{equation}\label{4_2}
    \phi_{\mathcal{B}} (A, X) := \inf \Big\{ \phi_1 (A, X, S^{*}) : \beta^{*} \in \mathcal{B} \Big\} > 0. 
\end{equation}
This quantity, with a special choice of $\mathcal{B}$ as laid out in the ensuing discussion, plays an important role in both Corollary 1 and Theorem \ref{thm:3}. We now turn our attention to the IC(Dimension) assumption. \newline
 
\noindent \textbf{IC (Dimension) Assumption:} There exists at least one $s \in \{ 3 b_n, \dots d_n \}$, and a corresponding quantity $\phi_0(A, X,s)>0$, such that for any $\beta_1, \beta_2 \in \mathbb{R}^{d_n}$, we have
$$\beta_{1} \neq \beta_{2}, |\supp{(\beta_2 - \beta_1)}| \leq s \quad \mathrm{implies} \quad \mathcal{D}_n(\eta_1| |\eta_2) \geq \frac{n \phi_0^2(A, X, s)}{|\supp{(\beta_2 - \beta_1)}|}. \| \beta_{2} - \beta_{1} \|_1^2.$$
The suffix $0$ of $\phi_0$ emphasizes we are working with constraints in the $\ell_0$ distance.  Similar to IC(Model), the intuition behind IC(Dimension) is to guarantee that whenever a general $\beta$ matches on most of the co-ordinates with true $\beta^{*}$, i.e. their $\ell_0$ distance is small, the KL term should be able to separate them. \newline

The IC(Dimension) assumption, coupled with the IC(Model) assumption, form one of the  central conditions in the proof of our posterior contraction statement, and we shall call it the IC (Joint) condition. First, consider the set
$$\mathcal{B}_{0, n} := \left\{ \beta \in \mathbb{R}^{d_n}: \overline{\phi}_0 \left(A, X, 3 |\supp{(\beta)}| \right) > 0 \right\},$$
where, for any $s \in \{ 1, \dots d_n \}$,
\begin{equation}\label{4_3}
    \overline{\phi}_0 (A, X, s) := \inf \left\{ \frac{\sqrt{ s \mathcal{D}_n(\eta^{*}| |\eta)}}{ \sqrt{n} \| \beta - \beta^{*} \|_1} : |\supp{(\beta- \beta^{*})}| \leq s, \; \beta \neq \beta^{*} \right\}.
\end{equation}
Now observe that $ \overline{\phi}_0(A, X, s)$ is decreasing in $s$, by definition, for any fixed $A$ and $X$. Now, by IC(Dimension), we have $\overline{\phi}_0(A, X, 3 b_n) > 0$, which shows $\overline{\phi}_0(A, X, 3|\supp{(\beta)}|) > 0$ whenever $|\supp{(\beta)}| \leq b_n$. We thus have
\begin{equation}\label{4_4}
    \mathcal{B}_{0, n} \supset \left\{ \beta \in \mathbb{R}^{d_n} : 0 < |\supp{(\beta)}| \leq b_n \right\} =: \mathcal{B}_{2, n},
\end{equation}
which is a desirable relation based on the discussion at the start of this section. We are now ready to state \newline

\noindent \textbf{IC (Joint) Assumption:}
$$\mathcal{B}_n := \mathcal{B}_{1, n} \cap \mathcal{B}_{2, n} \; \mathrm{is} \; \mbox{non-empty}.$$
A direct and vital consequence of this assumption is $\phi_{\mathcal{B}_n} (A, X) > 0$, as seen from \eqref{4_2} by choosing $\mathcal{B} = \mathcal{B}_n$. As we shall see, the statements of our results in Theorem \ref{thm:2} and Theorem \ref{thm:3} are uniformly adaptive over $\beta^{*} \in \mathcal{B}_n$. More precisely, our posterior contraction statement will have the form
$$\underset{\beta^{*} \in \mathcal{B}_n}{\sup} \bbP \left( \|\beta - \beta^{*}\|_1 > \varepsilon_{n, 1} \mid Y^{(n)} \right) \to 0 \quad \mathrm{as} \quad n \to \infty.$$
with $\varepsilon_{n, 1} > 0$ generically denoting the optimal radius of posterior contraction. \newline

We end this section with the pivotal role of clipping function $\eta(\cdot)$ in the IC assumptions. We require the geometries of the likelihood and the prior to match up in terms of the parameter $\beta$; $\mathcal{D}_n(\eta_1| |\eta_2)$ captures the discrepancy among $\beta$'s in the likelihood, while the $\ell_1$ gap does the same for the Laplace signals in the prior. To have a posterior contraction statement in $\ell_1$ distance, it is necessary for the $\mathcal{D}_n(\eta_1| |\eta_2)$ to grow with $\| \beta_2 - \beta_1 \|_1$, at least in sparsity restricted sense, and that is what the IC(Model) and IC(Dimension) assumptions reflect. Clipping function $\eta(\cdot)$, being an intermediary of $A(\cdot)$ and linear term $x_i^{\T} \beta$, must also reflect this growth, and hence has to be necessarily injective. The Lipschitz nature of $\eta(\cdot)$ allows us to translate gaps between $\eta$'s to gaps between $\beta$'s. \newline

\subsection{Hyper-parameter selection aimed at truth adaptive posterior contraction}
Since we aim to avoid prior dependence on the truth, choosing the hyper-parameter $\lambda_n, a_n$ should only take into account the sample size $n$, parameter dimension $d_n$, covariate matrix $X$ and log-partition function $A(\cdot)$. Our assumptions must allow us to forgo use of any prior knowledge of the truth $\beta^{*}$ while hyper-parameter selection. Choice of $\lambda_n$ is significantly inter-twined with the log-partition function $A(\cdot)$ as well as the clipping function $\eta(\cdot)$. As in \cite{castillo2015bayesian}, $\lambda_n$ needs to scale with some function of the design matrix $X$, and since covariate information from $X$ is fed into the log-partition function through $\eta(\cdot)$, choice of $\lambda_n$ depends on $A(\cdot), \eta(\cdot)$ and $X$. Based on this, consider the bound
$$\underset{\beta^{*} \in \mathcal{B}_{2,n}}{\sup} \underset{1 \leq i \leq n}{\max} \; \sup \left\{ A^{''} (\gamma) : \left|\gamma - \eta_i^{*} \right| \leq \sqrt{ \frac{s_n^{*} \log d_n}{n} } \right\} \leq \mathcal{M}_0^2(A),$$
which essentially gives us local control over $A''(\eta_i) \; \forall \; i = 1, \dots n$, uniformly over $\beta^{*} \in \mathbb{R}^{d_n}$. The proof of this statement is detailed in Lemma \ref{lemma:3} in the appendix, and basically uses two main points. Firstly, since $s_n^{*} \leq b_n$ for $\beta^{*} \in \mathcal{B}_n$, we have $(s_n^{*} \log d_n)/n \to 0$ by \eqref{4_1}, which allows us to have shrinking neighborhoods around every $\eta_i^{*}$. Secondly, based on the behavior of $A''(\cdot)$, the clipping function $\eta(\cdot)$ restricts the set of arguments passed to $A(\cdot)$, thus controlling the growth of $A''(\cdot)$. \newline

Now, define the quantities 
\begin{equation}\label{4_5}
    \begin{split}
        \mathcal{M}_1(A) &:= \left( 1 \wedge \mathcal{M}_0^{-1}(A) \right)^{-1}, \|X\|_{(\infty,\infty)} := \max \left\{ X_{i.j} : i=1, \dots n, j=1. \dots d_n \right\}, \\
        \mathcal{M}(A, X) &:= \|X\|_{(\infty,\infty)} \mathcal{M}_1(A).
    \end{split}
\end{equation}
We can now state our assumption on the hyper-parameter $\lambda_n$: \newline

\noindent \textbf{Assumption $\mathcal{L}_0$:}
$$\frac{\mathcal{M}(A, X)}{d_n} \leq \lambda_n \leq \mathcal{M}(A, X) \sqrt{\log d_n}.$$
This bound, which we utilize in all our Theorems, generalizes the hyper-parameter bounds mentioned in \cite{castillo2015bayesian}, as well as avoids prior dependence on the truth. Existence of $\mathcal{M}_0(A) > 0$, through which $\mathcal{M}(A, X)$ is defined in \eqref{4_5}, is guaranteed by \eqref{2_3}, and it acts as a pre-fixed constant quantity that the practitioner can choose based solely on $A(\cdot)$, and then choose clipping function $\eta(\cdot)$. This, in turn, shows that the choice of hyper-parameter $\lambda_n$ depends solely on the three quantities $(A, X, b_n)$. This makes our hyper-parameter choice of $\lambda_n$ free of the truth. \newline

We turn our attention to hyper-parameter $a_n$, which controls how fast the the model weights $\omega_n(\cdot)$ decay. First, define
\begin{equation}\label{4_6}
    \mathcal{E}_1 := 8 \left( 1 +  \frac{49 \mathcal{M}^2(A, X) }{8 \phi_{\mathcal{B}_n}^2 (A, X)} \right),
\end{equation}
which is an adaptive choice, as well as free of any knowledge of the truth, owing to \eqref{4_2} and IC (Joint) assumption. For mild demands, like in Theorem \ref{thm:2}, $a_n > 1$ suffices. On the contrary, for the weak model selection result in Corollary \ref{cor}, we need to choose $a_n$ that supports very strong down-weighting of larger models, namely $a_n \geq 1 + 2 b_n \mathcal{E}_1$. One can note from \eqref{4_1} why this choice of $a_n$ heaviliy penalizes larger models. Lastly, the choice $a_n \geq 1 + \mathcal{E}_1$, which is much milder than our previous choice, is sufficient for the posterior contraction result in Theorem \ref{thm:3}. It is crucial to note that just like $\lambda_n$, our choice of hyper-parameter $a_n$ avoids any knowledge of true $\beta^{*}$.
\section{Adaptive rate-optimal posterior contraction rate in $\ell_1$ norm}\label{sec:rate}
In this section, we provide the statements of our results and lay out sketches of how we arrive at them, putting under spotlight the use of the assumptions.
\subsection{Lower bound of the marginal likelihood}
Starting from \eqref{2_4}, the marginal likelihood is defined as
\begin{equation}\label{5_1}
     \int \exp \left( L_n(\eta, \eta^{*}) \right) \Pi_n(\beta) d\beta,
\end{equation}
which appears as the denominator in calculating the posterior through Bayes' Theorem. Theorem \ref{thm:1} provides a high probability lower bound to this quantity in terms of the parameter dimension $d_n$ and the true model size $s_n^{*}$. \newline

\begin{theorem}\label{thm:1} Let $a_n > 0$ and $\lambda_n$ satisfy assumption $\mathcal{L}_0$. Let $n, d_n \to \infty$ and $d_n > n$. Based on \eqref{4_1}, consider large enough $n$ so that $b_n \log d_n < n$. Let the true $\beta^{*}$ belong to $\mathcal{B}_{2, n}$ as in \eqref{4_4}. Then, with probability $1 - \left( s_n^{*} \log d_n \right)^{-1}$ with respect to the data generating distribution, the marginal likelihood defined in \eqref{5_1} satisfies for all sufficiently large $n$
$$\int \exp \left( L_n(\eta, \eta^{*}) \right) \Pi_n(\beta) d\beta  \gtrsim C_n d_n^{-(a_n+6)s_n^{*}} \exp (-\lambda_n \|\beta^{*}\|_1).$$
\end{theorem}
The fact that $s_n^{*} \log d_n \to \infty$ as $n \to \infty$ makes this a high probability statement about the marginal likelihood. Since majority of the mass under the integral should lie around the truth $\beta^{*}$, it is natural that the lower bound should contain information about that truth. Theorem \ref{thm:1} quantifies that relation. One generic tool for reaching such a bound has been described in \cite{ghosal2000convergence}, which we modify to suit our needs. We provide a small sketch of our method here, while the full proof is given in the appendix. \newline

Consider the set
\begin{equation}\label{5_2}
    D_n := \left\{ \beta \in \mathbb{R}^{d_n} : \Big[ - \E \left[ L_n(\eta, \eta^{*}) \right] \; \bigvee \; \Var \left[  L_n(\eta, \eta^{*}) \right] \Big] \leq s_n^{*} \log d_n \right\},
\end{equation}
noting that $ - \E \left[ L_n(\eta, \eta^{*}) \right] = \mathcal{D}_n (\eta^{*}| |\eta)$ and $\Var \left[  L_n(\eta, \eta^{*}) \right] = \E Z_n^2(\eta, \eta^{*})$. Let $\Pi_{D_n}(\beta)$ denote the restriction of prior $\Pi_n(\beta)$ to $D_n$. Then the denominator satisfies
$$\int \exp \left( L_n(\eta, \eta^{*}) \right) \Pi_n(\beta) d\beta \geq \int_{D_n} \exp \left( L_n(\eta, \eta^{*}) \right) \Pi_n(\beta) d\beta = \Pi_n \left( D_n \right) \int \exp \left( L_n(\eta, \eta^{*}) \right) \Pi_{D_n}(\beta) d\beta.$$
This method of restricting the integral of the marginal likelihood to a neighborhood of the truth is reminiscent of the original method found in \cite{ghosal2000convergence}.  The radius of such a neighborhood, here given by $s_n^{*} \log d_n$, signifies the order of allowable growth in both the expectation and variance of the log-likelihood ratio. We now have two terms to deal with, the prior probability of $D_n$ and the restricted integral. First, we use the variance of $L_n(\eta, \eta^{*})$ in a Chebyshev inequality to obtain the lower bound 
\begin{equation}\label{5_3}
    \int \exp \left( L_n(\eta, \eta^{*}) \right) \Pi_{D_n}(\beta) d\beta \geq d_n^{- 2 s_n^{*}}.
\end{equation}
with high probability. 
Next, as detailed in Lemma 1 in the appendix, we bound from below the prior probability of $D_n$ as
$$\Pi_n(D_n) \gtrsim C_n \exp (-\lambda_n \|\beta^{*}\|_1) d_n^{-(a_n+4)s_n^{*}}.$$
Theorem \ref{thm:1} now follows by combining Lemma \ref{lemma:1} with \eqref{5_3}.

\subsection{Posterior dimension and weak model selection}
We work with complexity priors that put increasingly higher penalty, or lower weight, on models that have larger sizes. It is expected that the posterior would reflect this prior property, which tantamounts to the posterior having vanishingly low probability of exceeding a certain dimension. Theorem \ref{thm:2} does exactly that, showing that the posterior should be at least as sparse as the true $\beta^{*}$, up to multiplicative constants. Sparsity is quantified using $|\supp{(\beta)}|$ and is compared with $s_n^{*}$, the true level of sparsity in $\beta^{*}$. \newline

\begin{theorem}\label{thm:2} Let $a_n > 1$ and $\lambda_n$ satisfy assumption $\mathcal{L}_0$. Let $n, d_n \to \infty$ and $d_n > n$. Based on \eqref{4_1}, consider large enough $n$ so that $b_n \log d_n < n$. Let assumptions IC(Model) and IC (Joint) hold, and consider the non-null set $\mathcal{B}_n$. Then, with quantity $\phi_1(A, X, S)$ given by IC(Model), and $\mathcal{M}(A, X)$ as in \eqref{4_5}, we have for all sufficiently large $n$,
$$\underset{\beta^{*} \in \mathcal{B}_n}{\sup} \E \Bigg[ \Pi_n \left( |\supp{(\beta)}| > s_n^{*} \left[1 + \frac{8}{a_n-1} \left( 1 +  \frac{49 \mathcal{M}^2(A, X)}{8 \phi_1^{2} (A, X, S^{*})} \right) \right] \Bigg| Y^{(n)} \right) \Bigg] \to 0 \quad \mathrm{as} \quad n \to \infty.$$
\end{theorem}
The statement of the theorem is presented in an asymptotic fashion, but is true for every $n$ large enough, satisfying the order assumptions. For simplicity, let us define the quantity $\mathcal{E}_1^{*} := 8 \left( 1 + 49 \mathcal{M}^2(A, X)/8 \phi_1^{2} (A, X, S^{*}) \right)$ so that Theorem \ref{thm:2} is a statement about the posterior probability of the set $\left\{ |\supp{(\beta)}| > s_n^{*} \left( 1 + \mathcal{E}_1^{*}/(a_n - 1) \right) \right\}$. It is important to note that we have used $\beta^{*} \in \mathcal{B}_n$ implies $\phi_1 (A, X, S^{*}) > 0$.  Owing to IC (Joint), \eqref{4_2} and the choice $\mathcal{B} = \mathcal{B}_n$, we can have from Theorem \ref{thm:2}, 
$$\underset{\beta^{*} \in \mathcal{B}_n}{\sup} \E \Bigg[ \Pi_n \left( |\supp{(\beta)}| > s_n^{*} \left[1 + \frac{8}{a_n-1} \left( 1 +  \frac{49 \mathcal{M}^2(A, X)}{8 \phi_{ \mathcal{B}_n}^2 (A, X)} \right) \right] \Bigg| Y^{(n)} \right) \Bigg] \to 0 \quad \mathrm{as} \quad n \to \infty.$$
By the definition of $\mathcal{E}_1$ in \eqref{4_6} and its analogy with $\mathcal{E}_1^{*}$, we now work with the posterior probability of $\left\{ |\supp{(\beta)}| > s_n^{*} \left( 1 + \mathcal{E}_1/(a_n - 1) \right) \right\}$. This allows to us to choose the hyper-parameter $a_n$ as $a_n \geq 1 + 2 b_n \mathcal{E}_1$, which is a truth-free choice, and leads to the following corollary: \newline

\begin{corollary}\label{cor}
With $\mathcal{E}_1$ as in \eqref{4_6}, if hyper-parameter $a_n$ in the prior satisfies $a_n \geq 1 + 2 b_n \mathcal{E}_1$ in addition to the hypotheses of Theorem \ref{thm:2}, we have
$$\underset{\beta^{*} \in \mathcal{B}_n}{\sup} \E \Bigg[ \Pi_n \left( \supp{(\beta)} \supsetneqq S^{*} \Big| Y^{(n)} \right) \Bigg] \to 0 \quad \mathrm{as} \quad n \to \infty.$$
\end{corollary}
This statement is a straightforward consequence of Theorem \ref{thm:2}, the fact that $s_n^{*} \leq b_n$ for $\beta^{*} \in \mathcal{B}_n$, and the observation that $ \left\{ \supp{(\beta)} \supsetneqq S^{*} \right\} \subset \left\{ |\supp{(\beta)}| > s_n^{*} + 1/2 \right\} $. Thus, Corollary 1 is a weak statement on model selection consistency. It ensures vanishingly small posterior probability attached to models that are strict super sets of the true model $S^{*}$. 

\subsection{Truth adaptive posterior contraction in $\ell_1$ metric}
We now turn our attention to the central result of our article, which is a truth adaptive statement about $\ell_1$-contraction of the posterior distribution. Essentially, it gives the radius of the smallest possible $\ell_1$ ball around true $\beta^{*}$, whose posterior probability vanishes with large $n$. Define the quantity
\begin{equation}\label{5_4}
    \mathcal{E}_2 := 6 + \frac{12 \mathcal{M}^2(A, X) }{ \overline{\phi}_0^2 \left(A, X, 3 b_n \right)},
\end{equation}
which can be observed to be truth-free. By describing the aforementioned radius in terms of $a_n, \mathcal{E}_2, d_n$ and $n$, we have the following.  
\begin{theorem}\label{thm:3}
Let hyper-parameter $a_n$ satisfy $a_n \geq 1 + \mathcal{E}_1$ for $\mathcal{E}_1$ as in \eqref{5_4}, and hyper-parameter $\lambda_n$ satisfy assumption $\mathcal{L}_0$. Let $n, d_n \to \infty$ and $d_n > n$. Based on \eqref{4_1}, consider large enough $n$ so that $b_n \log d_n < n$. Let assumptions IC(Model), IC(Dimension) and IC (Joint) hold, and consider the non-null set $\mathcal{B}_n$. Then, with quantity $\mathcal{E}_2$ given by \eqref{5_4} and $\mathcal{M}(A, X)$ as in \eqref{4_5}, we have  for all sufficiently large $n$,
\begin{equation}
    \underset{\beta^{*} \in \mathcal{B}_n}{\sup} \E \left[ \Pi_n \left( \| \beta - \beta^{*} \|_1 > \frac{2 s_n^{*} \left( 1 + a_n + \mathcal{E}_2 \right) }{\mathcal{M}(A, X)} \sqrt{\frac{\log d_n}{n}} \,\Big| \,  Y^{(n)} \right) \right] \to 0.
\end{equation}
\end{theorem}
It is important to note that the contraction rate linearly increases with $a_n$ and as long as $a_n$ is chosen to be a constant  larger than $1 + \mathcal{E}_1$, the rate is unaffected.  However, if one chooses a stronger penalty on the model space to achieve weak model selection consistency as in Corollary \ref{cor}, the rate of contraction in $\ell_1$ norm becomes slower unless the upper bound $b_n$ on the number of true non-zero coefficients  is assumed to be a constant. 

\section{Conclusion}
To summarize,  we introduced a new family of GLMs and developed sufficient conditions for obtaining posterior contraction rates that are adaptive rate-optimal.  From an implementation point of view, the new family does not bring in additional challenges, but from a theoretical point of view, it allows us to obtain adaptivity, while simultaneously obviating the need to enforce growth restriction on the true coefficient vector. Our analysis is restricted to the use of Laplace prior on the regression coefficients primarily due to the clarity and ease of calculations. More general priors, including compactly supported distributions, heavier tailed family or non-local priors can be considered. As a topic of immediate future research, strong model selection consistency is deemed important.  As already demonstrated in Theorem \ref{thm:2}, the posterior does not concentrate on subsets which are larger than the true subset with a stronger complexity prior.  With more identifiability conditions, one can ensure that the posterior does not concentrate on subsets that miss one or more non-zero true coordinates, thereby ensuring strong model selection consistency. 

\appendix
\section{Proof of Theorem \ref{thm:1}}\label{sec:thm1}
Recall the following neighborhood of the parameter space
$$D_n = \left\{ \beta \in \mathbb{R}^{d_n} : \Big[ - \E \left[ L_n(\eta, \eta^{*}) \right] \; \bigvee \; \Var \left[  L_n(\eta, \eta^{*}) \right] \Big] \leq s_n^{*} \log d_n \right\},$$
noting that $ - \E \left[ L_n(\eta, \eta^{*}) \right] = \mathcal{D}_n (\eta^{*}| |\eta)$ and $\Var \left[  L_n(\eta, \eta^{*}) \right] = \E Z_n^2(\eta, \eta^{*})$. Let $\Pi_{D_n}(\beta)$ denote the restriction of prior $\Pi_n(\beta)$ to $D_n$. Then the denominator satisfies
\begin{equation}\label{A1}
    \int \exp \left( L_n(\eta, \eta^{*}) \right) \Pi_n(\beta) d\beta \geq \int_{D_n} \exp \left( L_n(\eta, \eta^{*}) \right) \Pi_n(\beta) d\beta = \Pi_n \left( D_n \right) \int \exp \left( L_n(\eta, \eta^{*}) \right) \Pi_{D_n}(\beta) d\beta.
\end{equation}
We shall separately lower bound the the prior probability term and the integral term. First, we work with the integral in the above display. Rewrite $d_n^{-2s_n^{*}} = \exp \left( -2s_n^{*} \log d_n \right)$. Then consider following the tail event and inclusions, corresponding to the integral above:
\begin{equation*}
    \begin{split}
        \left\{ \int \exp \left( L_n(\eta, \eta^{*}) \right) \Pi_{D_n}(\beta) d\beta \leq d_n^{- 2 s_n^{*}} \right\}   \subset &\left\{  \int \Big[ Z_n (\eta, \eta^{*}) - \mathcal{D}_n (\eta^{*}| |\eta) \Big] \Pi_{D_n}(\beta) d\beta \leq -2 s_n^{*} \log d_n \right\} \\
        \subset &\left\{ \int Z_n (\eta, \eta^{*}) \Pi_{D_n}(\beta) \leq - s_n^{*} \log d_n \right\}.
    \end{split}
\end{equation*}
The first inclusion follows from Jensen's inequality, while the second uses that $\mathcal{D}_n (\eta^{*}| |\eta) \leq s_n^{*} \log d_n$ on $D_n$. We can now make the following probability statement about the integral,
\begin{equation}\label{A2}
    \begin{split}
        & \bbP \Big[ \int \exp \left( L_n(\eta, \eta^{*}) \right) \Pi_{D_n}(\beta) d\beta \leq d_n^{- 2 s_n^{*}} \Big] \leq \bbP \Big[ \int Z_n (\eta, \eta^{*}) \Pi_{D_n}(\beta) \leq - s_n^{*} \log d_n \Big] \\
        \leq &\frac{ \E \left[ \int Z_n(\eta, \eta^{*}) \Pi_{D_n}(\beta) d\beta \right]^2 }{\left( s_n^{*} \log d_n \right)^2} \leq \frac{\int \E Z_n^2(\eta, \eta^{*}) \Pi_{D_n}(\beta) d\beta}{\left( s_n^{*} \log d_n \right)^2} \leq \left( s_n^{*} \log d_n \right)^{-1},
    \end{split}
\end{equation}
where we have used Chebysev's inequality for the second inequality, the variance inequality in the third, and the fact that $ \E Z_n^2(\eta, \eta^{*}) \leq s_n^{*} \log d_n$ on $D_n$ for the final inequality. Note that \eqref{A2} makes sense asymptotically because $s_n^{*} \log d_n \to \infty$ with $n \to \infty$. Now recall the definition of $\mathcal{B}_{2, n}$ in \eqref{4_4}. Due to the hypothesis of Theorem \ref{thm:2}, we can use Lemma 1 to have
$$\Pi_n(D_n) \gtrsim C_n \exp (-\lambda_n \|\beta^{*}\|_1) d_n^{-(a_n + 4)s_n^{*}}.$$
Combining this with \eqref{A1} and \eqref{A2} shows that with probability greater than $1 - \left( s_n^{*} \log d_n \right)^{-1}$ w.r.t the data generating distribution for $\beta^{*} \in \mathcal{B}_{2, n}$, we have
$$\int \exp \left( L_n(\eta, \eta^{*}) \right) \Pi_n(\beta) d\beta \gtrsim C_n \exp (-\lambda_n \|\beta^{*}\|_1) d_n^{-(a_n + 6)s_n^{*}},$$
concluding the proof. 

\section{Proof of Theorem \ref{thm:2}}\label{sec:thm2}
We work with the quantity
$$\E \left[ \Pi_n(B \mid Y^{(n)}) \right] = \E \left[ \frac{\int_B \exp \left( L_n(\eta, \eta^{*}) \right) \Pi_n(\beta) d\beta}{\int \exp \left( L_n(\eta, \eta^{*}) \right) \Pi_n(\beta) d\beta} \right],$$
where the expectation is taken with respect to the true data generating distribution and set $B$ has the form $B = \left\{ \beta \in \mathbb{R}^{d_n} : | \supp{(\beta)} | > \varepsilon_n \right\}$ for some constant $\varepsilon_n > 0$. Thus, Theorem \ref{thm:2} is concerned with the dimensionality of the posterior vector, specifically the posterior probability that the sparsity of $\beta$ does not fall below a certain threshold. Start by defining $U_n := \mathcal{M}(A, X) \sqrt{n \log d_n}$, so that we have $\lambda_n \leq U_n$ from the hyper-parameter bounds assumption. Consider $\Omega_n:= \left\{ Y^{(n)}: Z_n(\eta, \eta^{*}) \leq U_n \| \beta - \beta^{*} \|_1 \right\}$, where $\Omega_n^c$ represents a tail event of the centered log-likelihood ratio $Z_n(\eta, \eta^{*})$. To find the probability of this event, observe that similar to calculations in Lemma 1,
$$\Var \left( Z_n(\eta, \eta^{*}) \right) = \sum_{i=1}^n \left( \eta_i - \eta_i^{*} \right)^2 A^{''} \left( \eta_i^{*} \right) \leq n \mathcal{M}^2(A, X)  \|\beta - \beta^{*} \|_1^2.$$
This shows, with the use of the definition of $U_n$ and Chebysev's inequality,
$$\bbP \left( \Omega_n^c \right) \leq \frac{\Var \left( Z_n(\eta, \eta^{*}) \right)}{U_n^2 \| \beta - \beta^{*} \|_1^2} \leq \left( \log d_n \right)^{-1}.$$
Also, Theorem \ref{thm:1} claims existence of event $\overline{\Omega}_n$ so that we have $ \bbP \left( \overline{\Omega}_n \right) \geq 1 - \left( s_n^{*} \log d_n \right)^{-1}$ and $\int \exp \left( L_n(\eta, \eta^{*}) \right) \Pi_n(\beta) d\beta  \gtrsim C_n d_n^{-(a_n + 6)s_n^{*}} \exp (-\lambda_n \|\beta^{*}\|_1)$ on $\overline{\Omega}_n$, simultaneously. Thus, using $\Pi_n(B \mid Y^{(n)}) \leq 1$ and the union bound for probabilities, we have
\begin{equation}\label{B1}
    \begin{split}
        &\hspace{12pt} \E \left[ \Pi_n(B \mid Y^{(n)}) \right] \leq \E \left[ \Pi_n(B \mid Y^{(n)}) \boldsymbol{1}_{\Omega_n \cap \overline{\Omega}_n} \right] + \bbP \left( \Omega_n^c \right) + \bbP \left( \overline{\Omega}_n^c \right) \\
        & \lesssim C_n^{-1} d_n^{(a_n + 6)s_n^{*}} \exp (\lambda_n \|\beta^{*}\|_1) \E \int_{B} \exp [L_n(\eta, \eta^{*})] \boldsymbol{1}_{\Omega_n} \Pi_n(\beta) d \beta + (\log d_n)^{-1} + \left( s_n^{*} \log d_n \right)^{-1},
    \end{split}
\end{equation}
since $\boldsymbol{1}_{\Omega_n \cap \overline{\Omega}_n} \leq \boldsymbol{1}_{\Omega_n} $. Since $\log d_n \to \infty$ by the order assumptions, it now suffices to work with the expectation term on the right hand side. Due to restriction to $\Omega_n$, we get
\begin{equation}\label{B2}
    \begin{split}
        & \E \int_{B} \exp [L_n(\eta, \eta^{*})] \boldsymbol{1}_{\Omega_n} \Pi_n(\beta) d \beta, \\
        &\leq \int_{B} \E \exp \left[ \left(1 - \frac{\lambda_n}{2 U_n} \right) Z_n(\eta, \eta^{*}) + \frac{\lambda_n}{2} \| \beta - \beta^{*} \|_1 - \mathcal{D}_n (\eta^{*}| |\eta) \right] \Pi_n(\beta) d \beta, \\
        &= \int_{B} \exp \left[ \frac{\lambda_n}{2} \| \beta - \beta^{*} \|_1 - \mathcal{D}_n (\eta^{*}| |\eta) \right]. \E \left( \exp \left[ \left(1 - \frac{\lambda_n}{2 U_n} \right) Z_n(\eta, \eta^{*}) \right] \right) \Pi_n(\beta) d \beta.
    \end{split}
\end{equation}
To calculate the expectation in the above display, we shall use Lemma \ref{lemma:2}, which is concerned with the connection of the KL divergence $\mathcal{D}_n (\eta^{*}| |\eta)$ with the cumulant generating function(cgf) of the centered log-likelihood ratio $Z_n(\eta, \eta^{*})$. We use $\alpha =  1 - \lambda_n/(2 U_n)$ in Lemma \ref{lemma:2}, obtaining
\begin{equation}\label{B3}
    \E \left( \exp \left[ \left(1 - \frac{\lambda_n}{2 U_n} \right) Z_n(\eta, \eta^{*}) \right] \right) \leq \left(1 - \frac{\lambda_n}{2 U_n} \right) \mathcal{D}_n (\eta^{*}| |\eta).
\end{equation}
The fact $0 < \lambda_n/(2 U_n) < 1$, implied by assumption $\mathcal{L}_0$, has been crucially used here. Combining \eqref{B2} and \eqref{B3}, we have for the expectation in \eqref{B1}
\begin{equation*}
    \begin{split}
        & \E \int_{B} \exp [L_n(\eta, \eta^{*})] \boldsymbol{1}_{\Omega_n} \Pi_n(\beta) d \beta \\
        &\leq \int_{B} \exp \left[ \frac{\lambda_n}{2} \| \beta - \beta^{*} \|_1 - \mathcal{D}_n (\eta^{*}| |\eta) \right].  \exp \left[ \left(1 - \frac{\lambda_n}{2 U_n} \right) \mathcal{D}_n (\eta^{*}| |\eta) \right] \Pi_n(\beta) d \beta \\
        &\leq\int_{B} \exp \left[ \frac{\lambda_n}{2} \| \beta - \beta^{*} \|_1 - \frac{\lambda_n}{2 U_n} \mathcal{D}_n (\eta^{*}| |\eta) \right] \Pi_n(\beta) d \beta,
    \end{split}
\end{equation*}
and hence
\begin{equation}\label{B4}
    \begin{split}
        &\hspace{12pt} \exp (\lambda_n \|\beta^{*}\|_1) . \E \int_{B} \exp [L_n(\eta, \eta^{*})] \boldsymbol{1}_{\Omega_n} \Pi_n(\beta) d \beta \\
        &\lesssim \int_{B} \exp \left[ \lambda_n \| \beta^{*} \|_1 + \frac{\lambda_n}{2} \| \beta - \beta^{*} \|_1 - \frac{\lambda_n}{2 U_n} \mathcal{D}_n (\eta^{*}| |\eta) \right] \Pi_n(\beta) d \beta.
    \end{split}
\end{equation}
We now work with the exponent inside the integrand in \eqref{B4}. First, $\| \beta^{*}\|_1 + (1/2) \| \beta - \beta^{*}\|_1 \leq \| \beta_{S^*} \|_1 + (3/2) \| \beta_{S^{*}} - \beta^{*} \|_1 + \frac{1}{2} \| \beta_{S^{*c}} \|_1$.  If $\| \beta_{S^{*c}} \|_1 \geq 7 \| \beta_{S^{*}} - \beta^{*} \|_1$, then 
\begin{equation}\label{B5}
   \| \beta_{S^*} \|_1 + \frac{3}{2} \| \beta_{s_n^{*}} - \beta^{*} \|_1 + \frac{1}{2} \| \beta_{S^{*c}} \|_1 \leq -\frac{1}{4} \| \beta - \beta^{*}\|_1 + \| \beta \|_1,
\end{equation}
and if $\| \beta_{S^{*c}} \|_1 < 7 \| \beta_{S^{*}} - \beta^{*} \|_1$, then we use the IC(Model) assumption to get
\begin{equation}\label{B6}
    \begin{split}
        &\| \beta^{*}\|_1 + \frac{1}{2} \| \beta - \beta^{*}\|_1 \leq \| \beta_{S^*} \|_1 + \frac{7}{2} \| \beta_{S^{*}} - \beta^{*} \|_1 - 2 \| \beta_{S^{*}} - \beta^{*} \|_1 + \frac{1}{2} \| \beta_{S^{*c}} \|_1 \\
        &\leq \frac{7}{2} \frac{\sqrt{\mathcal{D}_n (\eta^{*}| |\eta)s_n^{*}}}{\sqrt{n} \phi_1(A, X, S^{*})} - \frac{1}{4} \| \beta - \beta^{*}\|_1 + \| \beta \|_1 \\
        &\leq \frac{49 U_n s_n^{*}}{8 n \phi_1^{2} (A, X, S^{*})} + \frac{1}{2 U_n} \mathcal{D}_n (\eta^{*}| |\eta) - \frac{1}{4} \| \beta - \beta^{*}\|_1 + \| \beta \|_1.
    \end{split}
\end{equation}
The fact that $\beta^{*} \in \mathcal{B}_n$ implies $\phi_1 (A, X, S^{*}) > 0$ is crucially used above. Combining the above two exhaustive cases, namely \eqref{B5} and \eqref{B6}, we get for the integral in \eqref{B4}
\begin{equation}\label{B7}
    \begin{split}
        &\hspace{12pt} \int_{B} \exp \left[ \lambda_n \| \beta^{*} \|_1 + \frac{\lambda_n}{2} \| \beta - \beta^{*} \|_1 - \frac{\lambda_n}{2 U_n} \mathcal{D}_n (\eta^{*}| |\eta) \right] \Pi_n(\beta) d \beta \\
        &\leq \exp \left( \frac{49  U_n^2 s_n^{*}}{8 n \phi_1^{2} (A, X, S^{*})} \right) \int_{B} \exp \left[ - \frac{\lambda_n}{4} \| \beta - \beta^{*} \|_1 + \lambda_n \| \beta \|_1 \right] \Pi_n(\beta) d \beta,
    \end{split}
\end{equation}
where we have also used $\lambda_n \leq U_n$. Now consider the integral in \eqref{B7} for $B = \left\{ \beta \in \mathbb{R}^{d_n} : \supp{(\beta)} > \varepsilon_n \right\}$. Writing out the prior fully, we have
\begin{equation}\label{B8}
    \begin{split}
        &\hspace{10pt} \int_{| \supp{(\beta)} | > \varepsilon_n} \exp \left[ - \frac{\lambda_n}{4} \| \beta - \beta^{*} \|_1 + \lambda_n \| \beta \|_1 \right] \Pi_n(\beta) d \beta \\
        &= \sum_{|S| > \varepsilon_n} C_n {d_n \choose |S|}^{-1} \left( \frac{\lambda_n}{2 d_n^{a_n}} \right)^{|S|} \int \exp \left[ - \frac{\lambda_n}{4} \| \beta_S - \beta^{*} \|_1 \right] d \beta_S \\
        &\leq  C_n. \sum_{s = \varepsilon_n}^{d_n} \left(4 d_n^{-a_n} \right)^{s} \lesssim C_n \left(4 d_n^{-a_n} \right)^{\varepsilon_n} \leq C_n \left(d_n^{-(a_n - 1)} \right)^{\varepsilon_n},
    \end{split}
\end{equation}
where, for the first inequality, we have used $\| \beta_S - \beta^{*} \|_1 \geq \| \beta_S - \beta_S^{*} \|_1$ before performing the Laplace density integral, and have used $d_n \geq 4$. Now putting together the bounds in \eqref{B8}, \eqref{B7} and \eqref{B4}, and using $U_n = \mathcal{M}(A, X) \sqrt{n \log d_n}$, the first term in right-hand side of \eqref{B1} satisfies
$$ \E \left[ \Pi_n(B \mid Y^{(n)}) \boldsymbol{1}_{\Omega_n \cap \overline{\Omega}_n}   \right] \lesssim \exp \Bigg[ \log d_n. \left[ (a_n + 6)s_n^{*} - (a_n - 1) \varepsilon_n + \frac{49 \mathcal{M}^2(A, X) s_n^{*}}{8 \phi_1^{2} (A, X, S^{*})} \right] \Bigg] \to 0,$$
as soon as
$$(a_n + 6)s_n^{*} - (a_n - 1) \varepsilon_n + \frac{49 \mathcal{M}^2(A, X) s_n^{*}}{8 \phi_1^{2} (A, X, S^{*})} \leq - s_n^{*}, \, \mbox{equivalently}\; \varepsilon_n \geq s_n^{*} + \frac{8 s_n^{*}}{a_n - 1} \Bigg[1 +  \frac{49  \mathcal{M}^2(A, X)}{8 \phi_1^{2} (A, X, S^{*})} \Bigg],$$
as $s_n^{*} \log d_n \to \infty$ as $n \to \infty$. The proof is now completed by observing that the same lower bound on $\varepsilon_n$ works for every $\beta^{*} \in \mathcal{B}_n$.

\section{Proof of Theorem \ref{thm:3}}\label{sec:thm3}
Theorem \ref{thm:3} deals with the $\ell_1$-distance based posterior contraction of $\beta$ towards $\beta^{*}$, specifically the posterior probability of a set of the form $B_1 = \left\{ \beta \in \mathbb{R}^{d_n} : \|\beta - \beta^{*}\|_1 > \varepsilon_{n, 1} \right\}$. With $\mathcal{E}_1$ as in \eqref{4_6} and the choice $a_n \geq 1 + \mathcal{E}_1$, observe
\begin{equation}\label{C1}
    \left\{ \beta \in \mathbb{R}^{d_n} : |\supp{(\beta)}| \leq s_n^{*} \left( 1 + \frac{ \mathcal{E}_1}{a-1} , \right) \right\} \subset \left\{ \beta \in \mathbb{R}^{d_n} : |\supp{(\beta)}| \leq 2 s_n^{*} \right\} =: B_2.
\end{equation}
Now put $\overline{B} := B_1 \cap B_2$. As a result of Theorem \ref{thm:2}, we shall can focus only on $ \E \left[ \Pi_n( \overline{B} \mid Y^{(n)}) \right]$. Observe that $|\supp{(\beta - \beta^{*})}| \leq 3 s_n^{*} $ on $\overline{B}$ due to triangle inequality and \eqref{C1}. Now recall the definitions of $\Omega_n$ and $\overline{\Omega}_n$. Since $Z_n(\eta, \eta^{*}) \leq U_n \| \beta - \beta^{*} \|_1$ on $\Omega_n$, and $\lambda_n \leq 2 U_n$ implies $\lambda_n \| \beta^{*} \|_1 \leq \lambda_n \| \beta \|_1 + 2 U_n \| \beta - \beta^{*} \|_1$, we have
\begin{equation}\label{C2}
    \begin{split}
        &\hspace{10pt} \E \left[ \Pi_n(\overline{B} \mid Y) \boldsymbol{1}_{\Omega_n \cap \overline{\Omega}_n} \right] \lesssim C_n^{-1} d_n^{(a_n + 6)s_n^{*}} \exp (\lambda_n \|\beta^{*}\|_1) \E \int_{\overline{B}} \exp [L_n(\eta, \eta^{*})] \boldsymbol{1}_{\Omega_n} \Pi_n(\beta) d \beta \\
        &\leq C_n^{-1} d_n^{(a_n + 6)s_n^{*}} \int_{\overline{B}} \exp \left[ 4 U_n \| \beta - \beta^{*} \|_1 -  \mathcal{D}_n (\eta^{*}| |\eta) - U_n \| \beta - \beta^{*} \|_1 + \lambda_n \| \beta \|_1 \right] \Pi_n(\beta) d \beta,
    \end{split}
\end{equation}
where, similar to \eqref{B1} we already know $\bbP \left( \Omega_n \cap \overline{\Omega}_n \right)^c \to 0$ as $n \to \infty$. We shall now require the use IC(Dimension) assumption. Recall the definition of $\overline{\phi}_0(A, X, s)$ in \eqref{4_3}, which shows it is decreasing in $s$. Hence we have $\overline{\phi}_0 \left( A, X, |\supp{(\beta - \beta^{*})}| \right) \geq \overline{\phi}_0 \left(A, X, 3 s_n^{*} \right)$ whenever $\beta \in \overline{B}$, and $\overline{\phi}_0 \left(A, X, 3 s_n^{*} \right) > 0$ on account of $\beta^{*} \in \mathcal{B}_n$ and IC(Dimension) assumption. This leads to
\begin{equation}\label{C3}
    4 U_n \| \beta - \beta^{*} \|_1 \leq\frac{ 4 U_n \sqrt{ 3 s_n^{*} \mathcal{D}_n(\eta^{*}| |\eta)}}{ \sqrt{n} \overline{\phi}_0 \left(A, X, 3 s_n^{*} \right)} \leq \frac{12 U_n^2 s_n^{*}}{ n \overline{\phi}_0^2 \left(A, X, 3 s_n^{*} \right)} + \mathcal{D}_n(\eta^{*}| |\eta).
\end{equation}
Combining \eqref{C2} and \eqref{C3} with inequalities $\| \beta - \beta^{*} \|_1 > \varepsilon_{n, 1}$ for $\beta \in \overline{B}$, and $U_n \geq \frac{\lambda_n}{4} + \frac{U_n}{2}$, we have
\begin{equation*}
    \begin{split}
        &\hspace{10pt} \E \left[ \Pi_n(\overline{B} \mid Y) \boldsymbol{1}_{\Omega_n \cap \overline{\Omega}_n} \right] \\
        &\lesssim C_n^{-1} d_n^{(a_n + 6)s_n^{*}} \exp \left( \frac{12 U_n^2 s_n^{*} }{n \overline{\phi}_0^2 \left(A, X, 3 s_n^{*} \right) } - \frac{U_n \varepsilon_{n, 1}}{2}  \right) \int_{\overline{B}} \exp \left[ - \frac{\lambda_n}{4} \| \beta - \beta^{*} \|_1 + \lambda_n \| \beta \|_1 \right] \Pi_n(\beta) d \beta.
    \end{split}
\end{equation*}
Similar to calculations in \eqref{sec:thm2}, we note
$$\int_{\overline{B}} \exp \left[ - \frac{\lambda_n}{4} \| \beta - \beta^{*} \|_1 + \lambda_n \| \beta \|_1 \right] \Pi_n(\beta) d \beta \lesssim \frac{C_n}{1 - 4 d_n^{-a_n}} \lesssim C_n,$$
for large enough $d_n$, hence sufficiently large $n$, leading to
\begin{equation}\label{C4}
    \E \left[ \Pi_n(\overline{B} \mid Y) \boldsymbol{1}_{\Omega_n \cap \overline{\Omega}_n}  \right] \lesssim \exp \left[ \log d_n \left( (a_n + 6) s_n^{*} - \frac{\sqrt{n} \mathcal{M}(A, X) \varepsilon_{n, 1}}{2 \sqrt{\log d_n}} \right) + \frac{12 U_n^2 s_n^{*} }{n \overline{\phi}_0^2 \left(A, X, 3 s_n^{*} \right) } \right].
\end{equation}
Now recall the definition of $\mathcal{E}_2$ in \eqref{5_4}. Since $\beta^{*} \in \mathcal{B}_n$, we put to use that $s_n^{*} \leq b_n$ and that $\overline{\phi}_0 \left(A, X, s \right)$ is monotonically decreasing in $s$, to get from \eqref{C4}
$$ \E \left[ \Pi_n(\overline{B} \mid Y) \boldsymbol{1}_{\Omega_n \cap \overline{\Omega}_n} \right] \lesssim \exp \left[ \log d_n \left( s_n^{*} \left( a_n + \mathcal{E}_2 \right) - \frac{\sqrt{n}  \mathcal{M}(A, X) \varepsilon_{n, 1} }{2 \sqrt{\log d_n}} \right) \right] \to 0,$$
as soon as
$$s_n^{*} \left( a_n + \mathcal{E}_2 \right) - \frac{\sqrt{n} \mathcal{M}(A, X) \varepsilon_{n, 1} }{2 \sqrt{\log d_n}} \leq -s_n^{*} \Longleftrightarrow  > \varepsilon_{n, 1} \geq \frac{2 s_n^{*} \left( 1 + a_n + \mathcal{E}_2 \right) }{\mathcal{M}(A, X)} \sqrt{\frac{\log d_n}{n}},$$
as $s_n^{*} \log d_n \to \infty$ with $n \to \infty$. The proof is now completed by observing that the same lower bound on $\varepsilon_{n, 1}$ works for every $\beta^{*} \in \mathcal{B}_n$.

\section{Auxiliary results}\label{sec:aux}
In this section, we note down three Lemmata used in the proofs of our Theorems. Lemma \ref{lemma:1} deals with lower bounding $\Pi_n(D_n)$ with $D_n$ defined in \eqref{5_2}. Lemma \ref{lemma:2} upper bounds the cumulant generating function of $Z_n (\eta, \eta^{*})$ (defined in \eqref{2_4}), in terms of $\mathcal{D}_n (\eta^{*}| |\eta)$ (also defined in \eqref{2_4}). Lastly, Lemma \ref{lemma:3} deals with a local upper bound on $A''(\cdot)$ that is uniform over all $\beta^{*} \in \mathcal{B}_{2, n}$ (see \eqref{4_4} for definition).
\begin{lemma}\label{lemma:1} Let $a_n > 0$ and $\lambda_n$ satisfy assumption $\mathcal{L}_0$. Let $n, d_n \to \infty$ and $d_n > n$. Based on \eqref{4_1}, consider large enough $n$ so that $b_n \log d_n < n$. With $\mathcal{B}_{2, n}$ defined in \eqref{4_4}, let the true $\beta^{*}$ belong to $\mathcal{B}_{2, n}$. Then, for the set $D_n$ defined in \eqref{5_2}, we have for large enough $n$
$$\Pi_n(D_n) \geq C_n e^{-1/2} \exp (-\lambda_n \|\beta^{*}\|_1) d_n^{-(a_n+4)s_n^{*}}.$$
\end{lemma}
\noindent \textit{Proof:} Begin by defining
\begin{equation}\label{D1}
    \begin{split}
        B_n^{*}(A, X) &:= \mathcal{M}^{-1}(A, X) \sqrt{\frac{s_n^{*} \log d_n}{n}} , \\
        \Delta_n &:= \left\{ \beta \in \mathbb{R}^{d_n} : \|\beta - \beta^{*} \|_1 \leq B_n^{*}(A, X) \right\}.
    \end{split}
\end{equation}
Consider any $\beta \in \Delta_n$. We have 
$$\left| \eta \left( x_i^{T} \beta \right) - \eta \left( x_i^{T} \beta^{*} \right) \right| \leq \|X\|_{(\infty,\infty)} \| \beta - \beta^{*} \|_1 \leq \frac{1}{\mathcal{M}_1(A)} \sqrt{\frac{s_n^{*} \log d_n}{n}} \leq \sqrt{\frac{s_n^{*} \log d_n}{n}},$$
for all $i = 1, \dots n$, since  $\eta$ is a Lipschitz function, and $\mathcal{M}_1(A) \geq 1$ because of \eqref{4_5}. This shows, by Lemma \ref{lemma:3}, that for all $i = 1, \dots n$, we have $A^{''} \left( \gamma \right) \leq \mathcal{M}_0^2(A)$ whenever $\gamma$ lies between $\eta_i$ and $\eta_i^{*}$. Now note that
$$\Var \left[  L_n(\eta, \eta^{*}) \right] = \sum_{i=1}^n \left( \eta_i - \eta_i^{*} \right)^2 A^{''} \left( \eta_i^{*} \right), \; - \E \left[  L_n(\eta, \eta^{*}) \right] = \frac{1}{2} \sum_{i=1}^n \left( \eta_i - \eta_i^{*} \right)^2 A^{''} \left( \Tilde{\eta}_i \right),$$
where for all $i \in 1, \dots n$, we have $\Tilde{\eta}_i$ lying between $\eta_i$ and $\eta_i^{*}$. Thus we have for any $\beta \in \Delta_n$
\begin{equation}\label{D2}
    \begin{split}
        - \E \left[ L_n(\eta, \eta^{*}) \right] \; \bigvee \; \Var \left[  L_n(\eta, \eta^{*}) \right] &\leq \mathcal{M}_0^2(A) \sum_{i=1}^n \left( \eta_i - \eta_i^{*} \right)^2 \leq \mathcal{M}_1^2(A) \sum_{i=1}^n \left( \eta_i - \eta_i^{*} \right)^2 \\
        &\leq n \|X\|^2_{(\infty, \infty)} \mathcal{M}_1^2(A)  \|\beta - \beta^{*} \|_1^2 \leq s_n^{*} \log d_n.
    \end{split}
\end{equation}
Taken together, \eqref{D2} and \eqref{5_2} imply $\beta \in D_n$, which implies $\Delta_n \subset D_n$, and hence $\Pi_n \left( D_n \right) \geq \Pi_n \left( \Delta_n \right)$. Restricting prior $\Pi_n$ to the true model $S^{*}$, we see
\begin{equation*}
    \begin{split}
        \Pi_n(\Delta_n) &\geq C_n {d_n \choose s_n^{*}}^{-1} \left( \frac{\lambda_n}{2 d_n^{a_n}} \right)^{s_n^{*}} \int \exp (-\lambda_n \|\beta_{S^{*}}\|_1) \boldsymbol{1}_{ \left\{ \|\beta_{S^{*}} - \beta^{*} \|_1 \leq B_n^{*}(A, X) \right\} } d \beta_{S^{*}} \\
        &\geq C_n {d_n \choose s_n^{*}}^{-1} \left( \frac{\lambda_n}{2 d_n^{a_n}} \right)^{s_n^{*}} \exp (-\lambda_n \|\beta^{*}\|_1) \int \exp (-\lambda_n \|\chi_{S^{*}}\|_1) \boldsymbol{1}_{ \left\{ \| \chi_{S^{*}} \|_1 \leq B_n^{*}(A, X) \right\} } d \chi_{S^{*}},
    \end{split}
\end{equation*}
with the change of variable $\chi_{S^{*}} := \beta_{S^{*}} - \beta^{*}$, applying triangle inequality and noting that $\| \beta_{S^{*}}^{*} \|_1 = \| \beta^{*} \|_1$. To lower bound the above integral, we use it's analogy with Poisson process calculations. If we denote by $\mathcal{P}_j, j \geq 1$ independently and identically distributed exponential random variables with rate parameter $\lambda_n$, then the above integral is identical to calculating the probability of the event that at least $s_n^{*}$ many occurrences of the Poisson process $ \left\{ \sum_{j=1}^m \mathcal{P}_j \right\}_{m=1}^{\infty}$ happen before time $B_n^{*}(A, X)$. This leads us to
\begin{equation}\label{D3}
     \begin{split}
         \Pi_n(\Delta_n) &\geq C_n {d_n \choose s_n^{*}}^{-1} d_n^{-a_n s_n^{*}} \exp (-\lambda_n \|\beta^{*}\|_1) \exp \Big[ - \lambda_n B_n^{*}(A, X) \Big] \sum_{j=s_n^{*}}^{\infty} \frac{ \Big[ \lambda_n B_n^{*}(A, X) \Big]^{j}}{j!} \\
         &\geq C_n {d_n \choose s_n^{*}}^{-1} d_n^{-a_n s_n^{*}} \exp (-\lambda_n \|\beta^{*}\|_1) \exp \Big[ - \lambda_n B_n^{*}(A, X) \Big] \frac{ \Big[ \lambda_n B_n^{*}(A, X) \Big]^{s_n^{*}}}{s_n^{*}!} \\
         &\geq C_n d_n^{-(a_n + 1) s_n^{*}} \exp (-\lambda_n \|\beta^{*}\|_1) \exp \left[ - \frac{\lambda_n \epsilon_n}{\mathcal{M}(A, X)} \right] \left(  \frac{\lambda_n \epsilon_n}{\mathcal{M}(A, X)} \right)^{s_n^{*}}.
     \end{split}
\end{equation}
where $\epsilon_n := \sqrt{(s_n^{*} \log d_n)/n}$, so that $B_n^{*}(A, X) = \epsilon_n/\mathcal{M}(A, X)$ by \eqref{D1}, and we have used ${d_n \choose s_n^{*}} s_n^{*}! \leq d_n^{s_n^{*}}$ for the last inequality. Note that $\beta^{*} \in \mathcal{B}_n$ implies $s_n^{*} \leq b_n$, which, coupled with $b_n \log d_n < n$ implies $\epsilon_n < 1$. Based on assumption $\mathcal{L}_0$ and $d_n > n$, observe that
$$\frac{\lambda_n \epsilon_n}{\mathcal{M}(A, X)} \leq \frac{1}{2} \; \Rightarrow \; \exp \left[ - \frac{\lambda_n \epsilon_n}{\mathcal{M}(A, X)} \right] \left(  \frac{\lambda_n \epsilon_n}{\mathcal{M}(A, X)} \right)^{s_n^{*}} \geq e^{-\frac{1}{2}}. \left( \sqrt{\frac{s_n^{*} \log d_n}{n d_n^2}} \right)^{s_n^{*}} \geq e^{-\frac{1}{2}} d_n^{-\frac{3}{2} s_n^{*}},$$
and since $d_n \to \infty$ and $\epsilon_n < 1$, we have
$$\frac{\lambda_n \epsilon_n}{\mathcal{M}(A, X)} \geq \frac{1}{2} \; \Longrightarrow \; \exp \left[ - \frac{\lambda_n \epsilon_n}{\mathcal{M}(A, X)} \right] \left(  \frac{\lambda_n \epsilon_n}{\mathcal{M}(A, X)} \right)^{s_n^{*}} \geq 2^{-s_n^{*}} \exp \left[ -\sqrt{\log d_n} \right] \geq e^{-\frac{1}{2}} d_n^{-\frac{3}{2} s_n^{*}},$$
where the last inequality holds for large enough $d_n$, hence for large enough $n$. Plugging these back into the lower bound on $ \Pi_n(\Delta_n)$ in \eqref{D3}, we arrive at the statement of the Lemma.\newline

\begin{lemma}\label{lemma:2} Let the centered log-likelihood ratio $Z_n (\eta, \eta^{*})$ and the Kullback--Leibler divergence term $\mathcal{D}_n (\eta^{*}| |\eta)$ be defined as in \eqref{2_4}. Then, for any $\alpha \in (0, 1)$, the cumulant generating function $\psi(\alpha) := \log \E \left[ \exp \left( \alpha Z_n (\eta, \eta^{*}) \right) \right]$ of $Z_n(\eta, \eta^{*})$ satisfies
$$\psi(\alpha) \leq \alpha \mathcal{D}_n (\eta^{*}| |\eta).$$
\end{lemma}
\noindent \textit{Proof:} This Lemma is concerned with the connection of the KL divergence $\mathcal{D}_n (\eta^{*}| |\eta)$ in cGLM models with the cumulant generating function(cgf) of the centered log-likelihood ratio $Z_n(\eta, \eta^{*})$. Start by fixing $i = 1, \dots n$ and let $\alpha \in (0, 1)$. Put the cgf at $\alpha$ of $Z_i(\eta_i, \eta_i^{*})$ as $\psi_i(\alpha) := \log \E \left[ \exp \left( \alpha Z_i (\eta_i, \eta_i^{*}) \right) \right]$. Since $y_i$ is a draw from the exponential family, we know from standard properties that $\E T_i = A' \left( \eta_i^{*} \right)$ and for any $b \in \mathbb{R}$, $\log \E \left[ \exp(b T_i) \right] = A \left( \eta_i^{*} + b \right) - A \left( \eta_i^{*} \right)$. Hence, we have
\begin{equation}\label{D4}
    \begin{split}
        \psi_i(\alpha) &= \log \E \left( \exp \Big[ \alpha \left( \eta_i - \eta_i^{*} \right) \left(T_i - \E T_i \right) \Big] \right) = -\alpha \left( \eta_i - \eta_i^{*} \right) \E T_i + \log \E \left( \exp \Big[ \alpha \left( \eta_i - \eta_i^{*} \right) T_i \Big] \right) \\
        &= A \left( \eta^{*} + \alpha \left( \eta_i - \eta_i^{*} \right) \right) - A \left( \eta^{*} \right) - \alpha \left( \eta_i - \eta_i^{*} \right) A' \left( \eta_i^{*} \right) = \mathcal{D}_i \left( \eta_i^{*}| |\alpha \eta_i + (1-\alpha)\eta_i^{*} \right).
    \end{split}
\end{equation}
Now, since $\alpha \in (0, 1)$, we can use the convexity of KL divergence to obtain for every $i = 1, \dots n$
\begin{equation}\label{D5}
    \mathcal{D}_i \left( \eta_i^{*}| |\alpha \eta_i + (1-\alpha)\eta_i^{*} \right) \leq \alpha \mathcal{D}_i (\eta_i^{*}| |\eta_i).
\end{equation}
Since cgf of sum of independent random variables equals sum of their cgf's, we can sum over $i = 1, \dots n$ both the sides of \eqref{D5} and use \eqref{D4} for each term to obtain the statement of the Lemma. \newline

\begin{lemma}\label{lemma:3} Let $\mathcal{B}_{2, n}$ defined by \eqref{4_4}, while $\eta(\cdot)$ and $\mathcal{M}_0(A)$ are as in \eqref{2_3}. Then, for large enough $n$, we have
$$\underset{\beta^{*} \in \mathcal{B}_{2, n}}{\sup} \underset{1 \leq i \leq n}{\max} \; \sup \left\{ A^{''} (\gamma) : \left|\gamma - \eta_i^{*} \right| \leq \sqrt{ \frac{s_n^{*} \log d_n}{n} } \right\} \leq \mathcal{M}_0^2(A).$$
\end{lemma}
\noindent \textit{Proof:} Start with the simpler case, where $\mathcal{M}_0(A)$ can be chosen based on $A(\cdot)$, so that $\mathcal{I}_A \left( \mathcal{M}_0^2(A)/2 \right) = \mathbb{R}$. This results in $A''(\cdot)$ having the global upper bound $\mathcal{M}_0^2(A)/2$ on its support, and hence the above display holds trivially. Next, assume $\mathcal{I}_A(b)$ is a strict interval subset of $\mathbb{R}$ for any $b > 0$. For our proof, we shall only deal with interval form $\mathcal{I}_A \left( \mathcal{M}_0^2(A)/2 \right) = (-\infty, z_2], \; \mathcal{I}_A \left( \mathcal{M}_0^2(A) \right) = (-\infty, z_1]$, where $z_1, z_2 \in \mathbb{R}, \; z_1 > z_2$. All other form of intervals can be dealt with essentially the same technique we use. \newline

By the definition of clipping function, we have $\eta_i \in (- \infty, z_2], i = 1, \dots n$ for any $\beta \in \mathbb{R}^{d_n}$, specifically for any $\beta = \beta^{*} \in \mathcal{B}_{2, n}$.  Define the following neighborhood union
\begin{equation}\label{D6}
    \mathcal{N}_n := \bigcup_{i=1}^n \left\{ \gamma \in \mathbb{R} : \left|\gamma - \eta_i^{*} \right| \leq \sqrt{ \frac{b_n \log d_n}{n} } \right\}.
\end{equation}
\eqref{D6} deals with neighborhoods of $\eta_i^{*}$, where $\eta_i^{*} \equiv \eta \left( x_i^{\T} \beta^{*} \right), \beta^{*} \in \mathcal{B}_n, i = 1, \dots n$. By \eqref{4_1}, these neighborhoods shrink to zero with large $n$. Since $s_n^{*} \leq b_n$ is implied by $\beta^{*} \in \mathcal{B}_{2, n}$, and $z_1 > z_2$, we have for large enough $n$,
$$\mathcal{N}_n \subset \mathcal{I}_A \left( \mathcal{M}_0^2(A) \right),$$
implying 
\begin{equation}\label{D7}
    \underset{1 \leq i \leq n}{\max} \; \sup \left\{ A^{''} (\gamma) : \left|\gamma - \eta_i^{*} \right| \leq \sqrt{ \frac{s_n^{*} \log d_n}{n} } \right\} \leq \mathcal{M}_0^2(A).
\end{equation}
As \eqref{D7} holds for any $\beta^{*} \in \mathcal{B}_{2, n}$, this concludes the proof.

\bibliography{references}
\bibliographystyle{plainnat}


\end{document}